\documentclass[10pt]{amsart}
\usepackage{graphicx}
\usepackage{amssymb}
\usepackage{amsmath}
\newtheorem{theorem}{Theorem}[section]

\newtheorem{corollary}[theorem]{Corollary}

\newenvironment{pot1}{{\bf Proof of Theorem \ref{thm3.1}.}}{\hfill\fbox{}\par\vspace{.2cm}}
\newenvironment{pot2}{{\bf Proof of Theorem \ref{thm2}.}}{\hfill\fbox{}\par\vspace{.2cm}}
\newenvironment{pot3}{{\bf Proof of Theorem \ref{thm3}.}}{\hfill\fbox{}\par\vspace{.2cm}}
\newenvironment{poc1}{{\bf Proof of Corollary \ref{coro3.2}.}}{\hfill\fbox{}\par\vspace{.2cm}}

\numberwithin{equation}{section}

\def\charf {\mbox{{\text 1}\kern-.24em {\text l}}}

\def\bea{\begin{eqnarray*}}
\def\eea{\end{eqnarray*}}
\def\be{\begin{eqnarray}}
\def\ee{\end{eqnarray}}

\begin{document}

\title[Comparison geometry]{Comparison geometry for integral generalized quasi-Einstein tensor bounds}



\author{Sanghun Lee}
\address{Department of Mathematics, Chung-Ang University, 84 Heukseok-ro, Dongjak-gu, Seoul, Republic of Korea}
\curraddr{}
\email{kazauye@cau.ac.kr}


\subjclass[2010]{53C25; 53C20}

\keywords{Generalized quasi–Einstein tensor, Integral curvature, Comparison theorems, Diameter bounds}

\date{}

\dedicatory{}

\begin{abstract}
 The purpose of this paper is to extend the mean curvature comparison and volume comparison estimates by Petersen, Sprouse, and Wei to integral generalized quasi–Einstein tensor. Moreover, we use our comparison results to get global diameter estimate.
\end{abstract}

\maketitle
\section{Introduction}
 The study of generalized quasi–Einstein manifolds is introduced in \cite{CA}, which is a natural generalization of the Einstein metrics. More precisely, an $n$-dimensional complete Riemannian manifold $M$ is a \textit{generalized quasi–Einstein manifold} if there exist three smooth functions $f$, $\mu$, and $\lambda$ on $M$ satisfying
\be \label{eq1.1}
Ric_{f}^{\, \mu} = \lambda g.
\ee
Here,
\be \label{eq1.2}
Ric_{f}^{\, \mu} := Ric + \mbox{\rm Hess}f - \mu \, df\otimes df,
\ee 
where $Ric$ is the Ricci tensor on $M$ and \mbox{\rm Hess}$f$ is the Hessian of $f$. We call (\ref{eq1.2}) a \textit{generalized quasi–Einstein tensor}. When $\mu = \frac{1}{m}$ for a positive integer $m$, the above is called a generalized $m$-quasi–Einstein manifold. We denote $Ric_{f}^{\, m} := Ric + \mbox{\rm Hess}f - \frac{1}{m} \, df \otimes df$, and we call $Ric_{f}^{\, m}$ a \textit{generalized $m$-quasi–Einstein tensor}. The notion of generalized $m$-quasi–Einstein manifold was originated from the study of Einstein warped product manifolds (see \cite{BE}). It plays an important role in the study of the weighted measure (c.f. \cite{GW}). We also note that when $m = \infty$, a generalized $m$-quasi–Einstein tensor becomes
\be \label{eq1.3}
Ric_{f}^{\, \infty}:= Ric_{f} = Ric + \mbox{\rm Hess}f.
\ee
In particular, (\ref{eq1.3}) is called the \textit{Bakry–Emery Ricci tensor}.

 In \cite{PG1}, Petersen and Wei generalized the classical Bishop-Gromov volume comparison in an integral bound for the Ricci tensor. Before recalling their results, we need some notation. On a Riemannian manifold $M$, let $Ric_{-}$ be the smallest eigenvalue for the Ricci tensor $Ric: T_{x}M \rightarrow T_{x}M$
and
\bea
Ric^{\, H}_{-}:=\left((n-1)H - Ric_{-}\right)_{+} = \mbox{\rm max}\{0, (n-1)H - Ric_{-}\}
\eea
for $H \in \mathbb{R}$, the amount of Ricci tensor below $(n-1)H$. Also, we define 
\be \label{eq1.4.1}
\varphi := \left(m - m_{H}\right)_{+} = \mbox{\rm max}\{0, m - m_{H} \},
\ee
where $m$ is the mean curvature of the geodesic sphere in the outer normal direction and $m_{H}$ is the mean curvature of the geodesic sphere on the model space $M^{n}_{H}$. Here, $M^{n}_{H}$ is the $n$-dimensional simply connected space with constant sectional curvature $H$. We notice that if $r(x_{1}) = d(x_{1},x_{2})$ is the distance function from $x_{2}$ to $x_{1}$, then $\Delta(r) = m(r)$. Finally, we introduce a weighted $L^{p}$ norm of function $\phi$:
\bea
||\phi||_{p}(R) = \sup_{x \in M}\left(\int_{B(x,R)}|\phi|^{p} dv \right)^{\frac{1}{p}}.
\eea
Now we recall results of Petersen and Wei .

\begin{theorem} \cite{PG1}
 Let $M$ be an $n$-dimensional complete Riemannian manifold. For any $p > \frac{n}{2}$, $H \in \mathbb{R}$ (assume $r \leq \frac{\pi}{2\sqrt{H}}$ when $H > 0$),
\bea
||\varphi||_{2p}(r) \leq \left( \frac{(n-1)(2p-1)}{2p-n} ||Ric^{\,H}_{-}||_{p}(r) \right)^{\frac{1}{2}}.
\eea
Here, $\varphi$ is defined in (\ref{eq1.4.1}). Furthermore, for any $0<r<R$ (assume $R \leq \frac{\pi}{2\sqrt{H}}$ when $H>0$), there exists a constant $C(n,p,H,R)$ which is nondecreasing in $R$ such that
\bea
\left(\frac{V(x,R)}{V_{H}(R)}\right)^{\frac{1}{2p}} - \left(\frac{V(x,r)}{V_{H}(r)}\right)^{\frac{1}{2p}} \leq C(n,p,H,R)\left(||Ric^{\,H}_{-}||_{p}(R)||\right)^{\frac{1}{2}},
\eea 
where $V(x,R)$ is the volume of ball $B(x,R)$ in $M$, and $V_{H}(R)$ is the volume of ball $B(O,R)$ in the model space $M_{H}$.
\end{theorem}

 Using this comparison theorem, Petersen and Wei extended classical results such as Colding's volume convergence, Cheerger-Colding splitting theorem, Gromov precompactness theorem, and pinching result. Petersen and Sprouse extended the above comparison results and generalized Myers theorem in \cite{PC}. 

 On the other hand, Wu generalized Petersen, Sprouse, and Wei results \cite{PC, PG1, PG2} using Bakry–Emery Ricci tensor. In addition, Wu proved generalized Myers theorem, relative volume comparison for annulus, eigenvalue estimate, and volume growth estimate in \cite{JW1}. \\
 
 In this paper, we generalize Petersen, Sprouse, and Wei comparison results using to generalized quasi–Einstein tensor. Moreover, applying comparison results for generalized quasi–Einstein tensor, we prove diameter estimate. Given an $n$-dimensional smooth metric measure space $(M,g,e^{-f}dv)$, where $f$ is a smooth real valued function on $M$ and $dv$ is the usual Riemannian volume element on $M$. For the measure $e^{-f}dv$, the $f$-mean curvature is $m_{f} = m - \partial_{r}f$, where $m$ is the mean curvature of the geodesic sphere in the outer normal direction. The self-adjoint $f$-Laplacian with respect to the measure is $\Delta_{f} = \Delta - \nabla f\cdot\nabla$. Note that $m_{f}(r) = \Delta_{f}(r)$, where $r$ is the distance function.

 Let $Ric^{\, \mu}_{f-}$ be the smallest eigenvalue for the generalized quasi–Einstein tensor $Ric^{\, \mu}_{f}:T_{x}M \rightarrow T_{x}M$ and
\bea
Ric^{\, \mu, H}_{f-}:= \left((n+k-1)H - Ric^{\, \mu}_{f-}\right)_{+} = {\rm max} \, \{0, (n+k-1)H - Ric^{\, \mu}_{f-}\}
\eea
for $H \in \mathbb{R}$, the amount of $Ric^{\, \mu}_{f}$ below $(n-1)H$. We define a function $\varphi$ as follows:
\be \label{eq1.5.1}
\varphi := (m_{f} - m^{n+k}_{H})_{+} = {\rm max} \, \{0, m_{f} - m^{n+k}_{H} \},
\ee
where $m_{f} = m - \partial_{r}f$ and $m^{n+k}_{H}$ is the mean curvature of the geodesic sphere in the model space $M^{n+k}_{H}$. Also, we introduce a weighted $L^{p}$ norm of function $\phi$ on $(M,g,e^{-f}dv)$:
\bea
||\phi||_{p, \, f}(r) := \left( \int_{B(x,r)}|\phi|^{p}\mathcal{A}_{f}dtd\theta_{n-1}\right)^{\frac{1}{p}}.
\eea
Consider the quantity as
\bea
\bar{k}(p,H,r) := \sup_{x \in M}\left( \frac{1}{V_{f}(x,r)}\int_{B(x,r)}\left(Ric^{\, \mu,H}_{f-}\right)^{p}\mathcal{A}_{f}dtd\theta_{n-1}\right)^{\frac{1}{p}}.
\eea
Here, $\mathcal{A}_{f}(t,\theta)$ is the volume element of weighted form $e^{-f}dv = \mathcal{A}_{f}(t,\theta)dt\wedge d\theta_{n-1}$ in polar coordinate and $d\theta_{n-1}$ is the volume element on unit sphere $S^{n-1}$. Clearly, $Ric^{\, \mu}_{f} \geq (n+k-1)H$ if and only if $\bar{k}(p,H,r) = 0$. \\

 The following is a mean curvature comparison estimate for generalized quasi–Einstein tensor.

\begin{theorem} \label{thm1}
 Let $(M,g,e^{-f}dv)$ be an $n$-dimensional smooth metric measure space. Assume that $\mu \geq \frac{1}{k}$ for some positive constant $k$. For any $p>\frac{n+k}{2}$, $H \in \mathbb{R}$ (assume $r \leq \frac{\pi}{2\sqrt{H}}$ when $H>0$),
\begin{flalign} \label{eq1}
||\varphi||_{2p, \, f}(r) \leq \left(\frac{(n+k-1)(2p-1)}{2p-n-k}||Ric^{\, \mu, H}_{f-}||_{p, \, f}(r)\right)^{\frac{1}{2}}
\end{flalign}
and
\begin{flalign} \label{eq2}
\varphi^{2p-1}\mathcal{A}_{f} \leq (2p-1)^{p}\left(\frac{n+k-1}{2p-n-k}\right)^{p-1}\int^{r}_{0}\left(Ric^{\, \mu, H}_{f-}\right)^{p}\mathcal{A}_{f}\, dt.
\end{flalign}
Here, $\varphi$ is the function defined in (\ref{eq1.5.1}). Furthermore, if $H>0$ and $\frac{\pi}{2\sqrt{H}} < r < \frac{\pi}{\sqrt{H}}$, then we have
\begin{flalign} \label{eq1.6.3}
 ||sin^{\frac{4p-n-k-1}{2p}}(\sqrt{H}t)\varphi||_{2p, \, f}(r) \leq \left(\frac{(n+k-1)(2p-1)}{2p-n-k}||Ric^{\, \mu, H}_{f-}||_{p, \, f}(r)\right)^{\frac{1}{2}}
\end{flalign}
and
\begin{flalign} \label{eq1.6.4}
sin^{4p-n-k-1}(\sqrt{H}r)\varphi^{2p-1}\mathcal{A}_{f} \leq (2p-1)^{p}\left(\frac{n+k-1}{2p-n-k}\right)^{p-1}  \\
\times\int^{r}_{0}\left(Ric^{\, \mu, H}_{f-}\right)^{p}\mathcal{A}_{f}\, dt. \nonumber
\end{flalign}
\end{theorem}
 
 Secondly, we prove the following volume comparison estimate.

\begin{theorem} \label{thm2}
 Let $(M,g,e^{-f}dv)$ be an $n$-dimensional smooth metric measure space. Assume that $\mu \geq \frac{1}{k}$ for some positive constant $k$. For $H \in \mathbb{R}$, $p > \frac{n+k}{2}$, and $0 < r \leq R$ (assume $R \leq \frac{\pi}{2\sqrt{H}}$ when $H > 0$),
\bea
\left(\frac{V_{f}(x,R)}{V^{n+k}_{H}(R)}\right)^{\frac{1}{2p-1}} - \left(\frac{V_{f}(x,r)}{V^{n+k}_{H}(r)}\right)^{\frac{1}{2p-1}} &\leq& C(n+k,p,H,R) \\
&&\times \left(||Ric^{\, \mu,H}_{f-}||_{p, \, f}(R)\right)^{\frac{p}{2p-1}},
\eea
where 
\bea C(n+k,p,H,R) &:=& \left(\frac{n+k-1}{(2p-1)(2p-n-k)}\right)^{\frac{p-1}{2p-1}} \\
&&\times \displaystyle\int^{R}_{0}A^{n+k}_{H}(t)\left(\frac{t}{V^{n+k}_{H}(t)}\right)^{\frac{2p}{2p-1}}dt 
\eea
and $V^{n+k}_{H}(R)$ is the volume of ball $B(O,R)$ in the model space $M^{n+k}_{H}$ for $O \in M^{n+k}_{H}$.
\end{theorem}

 Using above theorems, we prove the following relative volume comparison for annulus as follows.

\begin{theorem} \label{thm3}
 Let $(M,g,e^{-f}dv)$ be an $n$-dimensional smooth metric measure space. Assume that $\mu \geq \frac{1}{k}$ for some positive constant $k$. Let $H \in \mathbb{R}$ and $p > \frac{n+k}{2}$. For $0 \leq r_{1} \leq r_{2} \leq R_{1} \leq R_{2}$ (assume $R_{2} \leq \frac{\pi}{2\sqrt{H}}$ when $H>0$),
\bea
\left(\frac{V_{f}(x,r_{2},R_{2})}{V^{n+k}_{H}(r_{2},R_{2})}\right)^{\frac{1}{2p-1}} - \left(\frac{V_{f}(x,r_{1},R_{1})}{V^{n+k}_{H}(r_{1},R_{1})}\right)^{\frac{1}{2p-1}} &\leq& C(n+k,p,H,r_{1},r_{2},R_{1},R_{2}) \\
&& \times \left(||Ric^{\, \mu, H}_{f-}||_{p, \, f}(R_{2})\right)^{\frac{p}{2p-1}},
\eea
where 
\begin{flalign*}
&C(n+k,p,H,r_{1},r_{2},R_{1},R_{2}) := \left(\frac{n+k-1}{(2p-1)(2p-n-k)}\right)^{\frac{p-1}{2p-1}}& \\
& \quad \times \left(\int^{R_{2}}_{R_{1}}A^{n+k}_{H}(t)\left(\frac{t}{V^{n+k}_{H}(r_{2},t)}\right)^{\frac{2p}{2p-1}}dt + \int^{r_{2}}_{r_{1}}A^{n+k}_{H}(R_{1})\left(\frac{R_{1}}{V^{n+k}_{H}(t,R_{1})}\right)^{\frac{2p}{2p-1}}dt\right).&
\end{flalign*}
\end{theorem}
 
 Finally, we apply the integral comparison results (Theorem \ref{thm1} and Theorem \ref{thm2}). That is, we obtain global diameter estimate for generalized quasi–Einstein tensor.

\begin{theorem} \label{thm4}
 Let $(M,g,e^{-f}dv)$ be an $n$-dimensional smooth metric measure space. Assume that $\mu \geq \frac{1}{k}$ for some positive constant $k$. Given $p>\frac{n+k}{2}$, $H > 0$, and $R>0$, there exist $D(n+k,H)$ and $\epsilon = \epsilon(n+k,p,H,R)$ such that if $\bar{k}(p,H,R) < \epsilon$, then $\mbox{\rm diam(M)} \leq D$.
\end{theorem}

 This paper is summarized as follows. In Section 2, we prove Theorem \ref{thm1}. In Section 3, we apply Theorem \ref{thm1} to prove Theorem \ref{thm2} and as corollary we prove volume doubling result when the integral generalized quasi–Einstein tensor bounds. Moreover, using Theorem \ref{thm2}, we obtain relative volume comparison for annulus. In Section 4, we give application of the integral comparison results. In other words, we apply Theorem \ref{thm1} and Theorem \ref{thm2} to get diameter estimate (Theorem \ref{thm4}).

\section{Mean curvaute comparison estimate}
 In this section, we prove mean curvature estimate for the integral generalized quasi–Einstein tensor. That is, we prove Theorem \ref{thm1}. The proof uses Bochner formula.
\\

 Applying the Bochner formula
\bea
 \frac{1}{2}\Delta|\nabla u|^{2} = |\mbox{\rm Hess} \, u|^{2} + \langle \nabla u, \nabla(\Delta u) \rangle + Ric(\nabla u, \nabla u)
\eea
to the distance function $r(x) = d(x,p)$, we obtain
\bea
0 = |\mbox{\rm Hess} \, r|^{2} + \frac{r}{\partial r}(\Delta r) + Ric(\nabla r, \nabla r).
\eea
Since $\mbox{\rm Hess} \,\, r$ is the second fundamental form of the geodesic sphere and $\Delta r$ is the mean curvaure, by the Schwarz inequality, we have
\bea
m' \leq -\frac{m^{2}}{n-1} - Ric(\nabla r, \nabla r).
\eea
And equality holds if and only if the radial sectional curvatures are constant. Hence the mean curvature of the $n$-dimensional model space $m_{H}$ satisfies
\bea
m'_{H} = -\frac{m^{2}_{H}}{n-1} - (n-1)H.
\eea
Since $m_{f} = m - \partial_{r}f$ and $m'_{f} = m' - \mbox{\rm Hess}\, f(\partial r, \partial r)$, we have
\bea
m'_{f} &\leq& -\frac{m^{2}}{n-1} - Ric(\partial r, \partial r) - \mbox{\rm Hess}\, f(\partial r, \partial r) \\
&=& -\frac{\left(m_{f} + \langle \nabla f, \nabla r \rangle (t)\right)^{2}}{n-1} - Ric(\partial r, \partial r) - \mbox{\rm Hess}\, f(\partial r, \partial r). 
\eea
Using the element inequality $(a+b)^{2} \geq \frac{a^{2}}{\alpha + 1} - \frac{b^{2}}{\alpha}$, for all real number $a$, $b$ and positive real number $\alpha$, we get
\bea
m'_{f} \leq - \frac{m^{2}_{f}}{(n-1)(\alpha + 1)} + \frac{\left(\langle \nabla f, \nabla r \rangle (t)\right)^{2}}{(n-1)\alpha} - Ric(\partial r, \partial r) - \mbox{\rm Hess}\, f(\partial r, \partial r).
\eea
Let $(n-1)\alpha = k$. Then we obtain
\bea
m'_{f} &\leq& - \frac{m^{2}_{f}}{n+k-1} + \frac{1}{k}\left(\langle \nabla f, \nabla r \rangle(t)\right)^{2} - Ric(\partial r, \partial r) - \mbox{\rm Hess}\, f(\partial r, \partial r) \\
&\leq& - \frac{m^{2}_{f}}{n+k-1} + \mu\left(\langle \nabla f, \nabla r \rangle(t)\right)^{2} - Ric(\partial r, \partial r) - \mbox{\rm Hess}\, f(\partial r, \partial r) \\ 
&=& - \frac{m^{2}_{f}}{n+k-1} - Ric^{\, \mu}_{f}(\partial r, \partial r).
\eea
We know that $\left(m^{n+k}_{H}\right)' = -\frac{\left(m^{n+k}_{H}\right)^{2}}{n+k-1} - (n+k-1)H$. So we compute
\bea
\left(m_{f} - m^{n+k}_{H}\right)' &\leq& - \frac{m^{2}_{f}}{n+k-1} - Ric^{\, \mu}_{f} + \frac{\left(m^{n+k}_{H}\right)^{2}}{n+k-1} + (n+k-1)H \\
&=& - \frac{\left(m_{f} - m^{n+k}_{H}\right)\left(m_{f} - m^{n+k}_{H} + 2m^{n+k}_{H}\right)}{n+k-1} \\
&& + (n+k-1)H - Ric^{\, \mu}_{f} \\
&\leq&  - \frac{\left(m_{f} - m^{n+k}_{H}\right)\left(m_{f} - m^{n+k}_{H} + 2m^{n+k}_{H}\right)}{n+k-1} + Ric^{\, \mu, H}_{f-}.
\eea
Note that on the interval where $m_{f} \leq m^{n+k}_{H}$, we have $\varphi = 0$ and where $m_{f} > m^{n+k}_{H}$, we have $\varphi = m_{f} - m^{n+k}_{H}$. Then we obtain
\bea
\varphi' + \frac{\varphi}{n+k-1}(\varphi + 2m^{n+k}_{H}) \leq Ric^{\, \mu, H}_{f-}.
\eea
Multiplying the above inequality by $(2p-1)\varphi^{2p-2}\mathcal{A}_{f}$, we get
\begin{flalign*}
&(2p-1)\varphi^{2p-2}\varphi'\mathcal{A}_{f} + \frac{2p-1}{n+k-1}\varphi^{2p}\mathcal{A}_{f} + \frac{4p-2}{n+k-1}\varphi^{2p-1}m^{n+k}_{H}\mathcal{A}_{f}& \\
&\quad \leq (2p-1)Ric^{\, \mu, H}_{f-}\varphi^{2p-2}\mathcal{A}_{f}.
\end{flalign*}
Using
\bea
(\varphi^{2p-1}\mathcal{A}_{f})' &=& (2p-1)\varphi^{2p-2}\varphi'\mathcal{A}_{f} + \varphi^{2p-1}\mathcal{A'}_{f} \\
&& (2p-1)\varphi^{2p-2}\varphi'\mathcal{A}_{f} + \varphi^{2p-1}m_{f}\mathcal{A}_{f},
\eea
we can be rewritten as
\begin{flalign*}
&(\varphi^{2p-1}\mathcal{A}_{f})' - \varphi^{2p-1}\mathcal{A}_{f}(m_{f} - m^{n+k}_{H}) - \varphi^{2p-1}\mathcal{A}_{f}m^{n+k}_{H} + \frac{2p-1}{n+k-1}\varphi^{2p}\mathcal{A}_{f}& \\
&\quad + \frac{4p-2}{n+k-1}\varphi^{2p-1}m^{n+k}_{H}\mathcal{A}_{f} \leq (2p-1)Ric^{\, \mu, H}_{f-}\varphi^{2p-2}\mathcal{A}_{f}.&
\end{flalign*}
We know that $- \varphi \leq -(m_{f} - m^{n+k}_{H})$. So we rearrange the above inequality as
\begin{flalign} \label{eq2.1}
&(\varphi^{2p-1}\mathcal{A}_{f})' + \varphi^{2p}\mathcal{A}_{f}\left(\frac{2p-1}{n+k-1} - 1 \right) + \varphi^{2p-1}m^{n+k}_{H}\mathcal{A}_{f}\left(\frac{4p-2}{n+k-1} - 1 \right)&   \\ 
&\quad \leq (2p-1)Ric^{\, \mu, H}_{f-}\varphi^{2p-2}\mathcal{A}_{f}.& \nonumber
\end{flalign}
Since $p > \frac{n+k}{2}$ and the assumption $r \leq \frac{\pi}{2\sqrt{H}}$, we have 
\bea
\varphi^{2p-1}m^{n+k}_{H}\mathcal{A}_{f}\left(\frac{4p-2}{n+k-1} - 1 \right) \geq 0.
\eea
Then we obtain
\bea
(\varphi^{2p-1}\mathcal{A}_{f})' + \varphi^{2p}\mathcal{A}_{f}\left(\frac{2p-n-k}{n+k-1} \right) \leq (2p-1)Ric^{\, \mu, H}_{f-}\varphi^{2p-2}\mathcal{A}_{f}.
\eea
Integrating the above inequality from $0$ to $r$, we get
\bea
\varphi^{2p-1}\mathcal{A}_{f} + \frac{2p-n-k}{n+k-1}\int^{r}_{0}\varphi^{2p}\mathcal{A}_{f} \, dt \leq (2p-1)\int^{r}_{0}Ric^{\, \mu, H}_{f-}\varphi^{2p-2}\mathcal{A}_{f} \, dt.
\eea
This implies that
\bea
\varphi^{2p-1}\mathcal{A}_{f} \leq (2p-1)\int^{r}_{0}Ric^{\, \mu, H}_{f-}\varphi^{2p-2}\mathcal{A}_{f} \, dt
\eea
and
\bea
\frac{2p-n-k}{n+k-1}\int^{r}_{0}\varphi^{2p}\mathcal{A}_{f} \, dt \leq (2p-1)\int^{r}_{0}Ric^{\, \mu, H}_{f-}\varphi^{2p-2}\mathcal{A}_{f} \, dt.
\eea
First, by Holder inequality, we have
\bea
\frac{2p-n-k}{n+k-1}\int^{r}_{0}\varphi^{2p}\mathcal{A}_{f} \, dt \leq (2p-1)\left(\int^{r}_{0}\left(Ric^{\, \mu, H}_{f-}\right)^{p}\mathcal{A}_{f} \, dt \right)^{\frac{1}{p}}\left(\int^{r}_{0}\varphi^{2p}\mathcal{A}_{f} \, dt \right)^{1 - \frac{1}{p}}.
\eea
Multiplying the above inequailty by $\frac{n+k-1}{2p-n-k}\left(\int^{r}_{0}\varphi^{2p}\mathcal{A}_{f} \, dt\right)^{-1 + \frac{1}{p}}$, we obtain
\begin{flalign} \label{eq2.2}
\left(\int^{r}_{0}\varphi^{2p}\mathcal{A}_{f} \, dt\right)^{\frac{1}{p}} \leq \frac{(n+k-1)(2p-1)}{2p-n-k}\left(\int^{r}_{0}\left(Ric^{\, \mu, H}_{f-}\right)^{p}\mathcal{A}_{f}\, dt\right)^{\frac{1}{p}}.
\end{flalign}
Take a root on both sides. Then we get (\ref{eq1}). \\
Secondly, we also have
\begin{flalign} \label{eq2.3}
\varphi^{2p-1}\mathcal{A}_{f} \leq (2p-1)\left(\int^{r}_{0}\left(Ric^{\, \mu, H}_{f-}\right)^{p}\mathcal{A}_{f} \, dt\right)^{\frac{1}{p}}\left(\int^{r}_{0}\varphi^{2p}\mathcal{A}_{f} \, dt\right)^{1 - \frac{1}{p}}.
\end{flalign}
Combining (\ref{eq2.2}) and (\ref{eq2.3}), we immediately yield (\ref{eq2}).

 If $H>0$ and $\frac{\pi}{2\sqrt{H}} < r < \frac{\pi}{\sqrt{H}}$, then $m^{n+k}_{H} < 0$. We know that
\begin{flalign*}
&(\varphi^{2p-1}\mathcal{A}_{f})' + \varphi^{2p}\mathcal{A}_{f}\left(\frac{2p-1}{n+k-1} - 1 \right) + \varphi^{2p-1}m^{n+k}_{H}\mathcal{A}_{f}\left(\frac{4p-2}{n+k-1} - 1 \right)& \\
&\quad \leq (2p-1)Ric^{\, \mu, H}_{f-}\varphi^{2p-2}\mathcal{A}_{f}.&
\end{flalign*}
Multiplying the above inequality by $sin^{4p-n-k-1}(\sqrt{H}r)$ and integrating from $0$ to $r$, we obtain
\begin{flalign*}
&\int^{r}_{0}sin^{4p-n-k-1}(\sqrt{H}t)\left(\varphi^{2p-1}\mathcal{A}_{f}\right)' dt + \frac{2p-n-k}{n+k-1}\int^{r}_{0}sin^{4p-n-k-1}(\sqrt{H}t)\varphi^{2p}\mathcal{A}_{f}\, dt& \\
&\quad + \frac{4p-n-k-1}{n+k-1}\int^{r}_{0}sin^{4p-n-k-1}(\sqrt{H}t)\varphi^{2p-1}m^{n+k}_{H}\mathcal{A}_{f}\, dt \\
&\qquad \leq (2p-1)\int^{r}_{0}sin^{4p-n-k-1}(\sqrt{H}t)Ric^{\, \mu, H}_{f-}\varphi^{2p-2}\mathcal{A}_{f} \, dt.
\end{flalign*}
Note that
\begin{flalign*}
&\int^{r}_{0}sin^{4p-n-k-1}(\sqrt{H}t)\left(\varphi^{2p-1}\mathcal{A}_{f}\right)' dt& \\
&\quad = - \frac{4p-n-k-1}{n+k-1}\int^{r}_{0}sin^{4p-n-k-1}(\sqrt{H}t)m^{n+k}_{H}\varphi^{2p-1}\mathcal{A}_{f} \, dt \\
&\qquad + sin^{4p-n-k-1}(\sqrt{H}r)\varphi^{2p-1}\mathcal{A}_{f}.
\end{flalign*}
So we have
\begin{flalign*}
&sin^{4p-n-k-1}(\sqrt{H}r)\varphi^{2p-1}\mathcal{A}_{f} + \frac{2p-n-k}{n+k-1}\int^{r}_{0}sin^{4p-n-k-1}(\sqrt{H}t)\varphi^{2p}\mathcal{A}_{f}\, dt& \\
&\quad \leq (2p-1)\int^{r}_{0}sin^{4p-n-k-1}(\sqrt{H}t)Ric^{\, \mu, H}_{f-}\varphi^{2p-2}\mathcal{A}_{f} \, dt.&
\end{flalign*}
This implies that
\begin{flalign} \label{eq2.4}
&sin^{4p-n-k-1}(\sqrt{H}r)\varphi^{2p-1}\mathcal{A}_{f}& \\
&\quad \leq (2p-1) \int^{r}_{0}sin^{4p-n-k-1}(\sqrt{H}t)Ric^{\, \mu, H}_{f-}\varphi^{2p-2}\mathcal{A}_{f} \, dt& \nonumber
\end{flalign}
and
\begin{flalign} \label{eq2.5}
&\frac{2p-n-k}{n+k-1}\int^{r}_{0}sin^{4p-n-k-1}(\sqrt{H}t)\varphi^{2p}\mathcal{A}_{f}\, dt& \\
&\quad \leq (2p-1)\int^{r}_{0}sin^{4p-n-k-1}(\sqrt{H}t)Ric^{\, \mu, H}_{f-}\varphi^{2p-2}\mathcal{A}_{f} \, dt.& \nonumber
\end{flalign}
Similar to the above discussion, by Holder inequality, we get
\begin{flalign} \label{eq2.6}
&\int^{r}_{0}sin^{4p-n-k-1}(\sqrt{H}t)Ric^{\, \mu, H}_{f-}\varphi^{2p-2}\mathcal{A}_{f} \, dt&  \\ 
&\quad \leq \left(\int^{r}_{0}sin^{4p-n-k-1}(\sqrt{H}t)\varphi^{2p}\mathcal{A}_{f} \, dt\right)^{1 - \frac{1}{p}}\nonumber \\
&\qquad \times \left(\int^{r}_{0}sin^{4p-n-k-1}(\sqrt{H}t)\left(Ric^{\, \mu, H}_{f-}\right)^{p}\mathcal{A}_{f} \, dt\right)^{\frac{1}{p}}. \nonumber
\end{flalign}
Substituting (\ref{eq2.6}) into (\ref{eq2.5}), then we obtain (\ref{eq1.6.3}). Similarly, putting (\ref{eq1.6.3}) and (\ref{eq2.6}) to (\ref{eq2.4}) immediately yields (\ref{eq1.6.4}).
\hfill\fbox{}\par\vspace{.5cm}

\section{Volume comparison estimates}
 In this section, we prove volume comparison results using mean curvature estimate. Before proving Theorem \ref{thm2}, we prove volume element comparison estimate when the integral generalized quasi–Einstein tensor bounds.

 For an $n$-dimensional smooth metric measure space $(M,g,e^{-f}dv)$, $\mathcal{A}_{f}(t,\theta)$ is the volume element of the weighted volume form $e^{-f}dv = \mathcal{A}_{f}(t,\theta)dt \wedge d\theta$ in the polar coordinate. In other words, $\mathcal{A}_{f}(t,\theta) = e^{-f}\mathcal{A}(t, \theta)$, where $\mathcal{A}(t,\theta)$ is the standard volume element of the metric $g$. Let
\bea
A_{f}(x,r) = \int_{S^{n-1}}\mathcal{A}_{f}(r,\theta)d\theta_{n-1}.
\eea
Similarly, we also let
\bea
A^{n+k}_{H}(r) = \int_{S^{n-1}}\mathcal{A}^{n+k}_{H}(r,\theta)d\theta,
\eea
where $\mathcal{A}^{n+k}_{H}$ is the volume element in the model space $M^{n+k}_{H}$.

 Now we prove volume element comparison estimate using the mean curvature estimate in Section 2.

\begin{theorem} \label{thm3.1}
 Let $(M, g, e^{-f}dv)$ be an $n$-dimensional smooth metric measure space. Assume that $\mu \geq \frac{1}{k}$ for some positive constant $k$. For $H \in \mathbb{R}$, $p > \frac{n+k}{2}$, and $0 < r \leq R$ (assume $R \leq \frac{\pi}{2\sqrt{H}}$ when $H >0$), we have
\begin{flalign} \label{eq3.1}
\left(\frac{A_{f}(x,R)}{A^{n+k}_{H}(R)}\right)^{\frac{1}{2p-1}} - \left(\frac{A_{f}(x,r)}{A^{n+k}_{H}(r)}\right)^{\frac{1}{2p-1}} &\leq \mathcal{C}(n+k,p,H,R)  \\
&\quad \times \left(||Ric^{\, \mu, H}_{f-}||_{p, \, f}(R)\right)^{\frac{p}{2p-1}}, \nonumber
\end{flalign}
where $\mathcal{C}(n+k,p,H,R) := \left(\frac{n+k-1}{(2p-n-k)(2p-1)}\right)^{\frac{p-1}{2p-1}}\int^{R}_{0}\left(A^{n+k}_{H}(t)\right)^{-\frac{1}{2p-1}}dt$. \\
Moreover, if $H>0$ and $\frac{\pi}{2\sqrt{H}} < r \leq R < \frac{\pi}{\sqrt{H}}$, then we obtain
\begin{flalign} \label{eq3.2}
&\left(\frac{A_{f}(x,R)}{A^{n+k}_{H}(R)}\right)^{\frac{1}{2p-1}} - \left(\frac{A_{f}(x,r)}{A^{n+k}_{H}(r)}\right)^{\frac{1}{2p-1}}& \\
&\quad \leq \left(\frac{n+k-1}{(2p-1)(2p-n-k)}\right)^{\frac{p-1}{2p-1}}\left(||Ric^{\, \mu, H}_{f-}||_{p, \, f}(R)\right)^{\frac{p}{2p-1}} \nonumber \\
&\qquad \times \int^{R}_{r}\frac{(\sqrt{H})^{\frac{n+k-1}{2p-1}}}{sin^{2}(t\sqrt{H})} \, dt. \nonumber
\end{flalign}
\end{theorem}

\begin{pot1}
 We apply $\mathcal{A'}_{f} = m_{f}\mathcal{A}_{f}$ and $(\mathcal{A}^{n+k}_{H})' = m^{n+k}_{H}\mathcal{A}^{n+k}_{H}$  to compute that
\bea
\frac{d}{dt}\left(\frac{\mathcal{A}_{f}(t,\theta)}{\mathcal{A}^{n+k}_{H}(t)}\right) &=& \frac{\mathcal{A}_{f}(t,\theta)}{\mathcal{A}^{n+k}_{H}(t)}(m_{f} - m^{n+k}_{H}) \\
&\leq&\frac{\mathcal{A}_{f}(t,\theta)}{\mathcal{A}^{n+k}_{H}(t)}\varphi.
\eea
So we have
\be \label{eq3.3}
\frac{d}{dt}\left(\frac{A_{f}(t,x)}{A^{n+k}_{H}(t)}\right) &=& \frac{1}{vol(S^{n-1})}\int_{S^{n-1}}\frac{d}{dt}\left(\frac{\mathcal{A}_{f}(t,\theta)}{\mathcal{A}^{n+k}_{H}(t)}\right)d\theta_{n-1}  \\
&\leq& \frac{1}{A^{n+k}_{H}(t)}\int_{S^{n-1}}\varphi \mathcal{A}_{f}(t,\theta) \, d\theta_{n-1}. \nonumber
\ee
By Holder inequality and (\ref{eq2}), we obtain
\begin{flalign} \label{eq3.4.1}
\frac{d}{dt}\left(\frac{A_{f}(t,x)}{A^{n+k}_{H}(t)}\right) \leq& \left((2p-1)^{p}\left(\frac{n+k-1}{2p-n-k}\right)^{p-1}\right)^{\frac{1}{2p-1}} \left(\frac{A_{f}(t,x)}{A^{n+k}_{H}(t)}\right)^{1 - \frac{1}{2p-1}} \\
&\times \left(||Ric^{\, \mu, H}_{f-}||_{p, \, f}(t)\right)^{\frac{p}{2p-1}}\left(\frac{1}{A^{n+k}_{H}(t)}\right)^{\frac{1}{2p-1}}. \nonumber
\end{flalign}
Note that
\bea
\frac{d}{dt}\left(\frac{A_{f}(t,x)}{A^{n+k}_{H}(t)}\right)\cdot\left(\frac{A_{f}(t,x)}{A^{n+k}_{H}(t)}\right)^{-1 + \frac{1}{2p-1}} = (2p-1)\frac{d}{dt}\left(\frac{A_{f}(t,x)}{A^{n+k}_{H}(t)}\right)^{\frac{1}{2p-1}}.
\eea
Then we get
\bea
\frac{d}{dt}\left(\frac{A_{f}(t,x)}{A^{n+k}_{H}(t)}\right)^{\frac{1}{2p-1}} &\leq& \left(\frac{n+k-1}{(2p-n-k)(2p-1)}\right)^{\frac{p-1}{2p-1}} \\
&&\times\left(||Ric^{\, \mu, H}_{f-}||_{p, \, f}(t)\right)^{\frac{p}{2p-1}}\left(\frac{1}{A^{n+k}_{H}(t)}\right)^{\frac{1}{2p-1}}.
\eea
Integrating from $r$ to $R$, we have
\begin{flalign*}
&\left(\frac{A_{f}(R,x)}{A^{n+k}_{H}(R)}\right)^{\frac{1}{2p-1}} - \left(\frac{A_{f}(r,x)}{A^{n+k}_{H}(r)}\right)^{\frac{1}{2p-1}}&\\
&\quad \leq \left(\frac{n+k-1}{(2p-n-k)(2p-1)}\right)^{\frac{p-1}{2p-1}}\left(||Ric^{\, \mu,H}_{f-}||_{p, \, f}(R)\right)^{\frac{p}{2p-1}}\int^{R}_{0} \left(\frac{1}{A^{n+k}_{H}(t)}\right)^{\frac{1}{2p-1}}dt.&
\end{flalign*}
For $H >0$ and $\frac{\pi}{2\sqrt{H}} < r \leq R < \frac{\pi}{\sqrt{H}}$, by (\ref{eq1.6.4}), we obtain
\begin{flalign} \label{eq3.4}
&\left(\int_{S^{n-1}}\varphi^{2p-1}\mathcal{A}_{f}d\theta_{n-1}\right)^{\frac{1}{2p-1}}& \\
&\quad \leq \left((2p-1)^{p}\left(\frac{n+k-1}{2p-n-k}\right)^{p-1}\right)^{\frac{1}{2p-1}}sin^{\frac{-4p+n+k+1}{2p-1}}(\sqrt{H}t) \nonumber  \\
&\qquad \times \left(||Ric^{\, \mu, H}_{f-}||_{p, \, f}(t)\right)^{\frac{p}{2p-1}}. \nonumber
\end{flalign}
Note that
\begin{flalign} \label{eq3.5}
\int_{S^{n-1}}\varphi \mathcal{A}_{f}(t,\theta) \, d\theta_{n-1} \leq \left(\int_{S^{n-1}}\varphi^{2p-1}\mathcal{A}_{f} \, d\theta_{n-1}\right)^{\frac{1}{2p-1}}A_{f}(t,x)^{1 - \frac{1}{2p-1}}.
\end{flalign}
Putting (\ref{eq3.4}) and (\ref{eq3.5}) into (\ref{eq3.3}), we get
\bea
\frac{d}{dt}\left(\frac{A_{f}(t,x)}{A^{n+k}_{H}(t)}\right) &\leq& \left((2p-1)^{p}\left(\frac{n+k-1}{2p-n-k}\right)^{p-1}\right)^{\frac{1}{2p-1}}sin^{\frac{-4p+n+k+1}{2p-1}}(\sqrt{H}t) \\
&& \times \left(||Ric^{\, \mu, H}_{f-}||_{p, \, f}(t)\right)^{\frac{p}{2p-1}}\left(\frac{A_{f}(t,x)}{A^{n+k}_{H}(t)}\right)^{1 - \frac{1}{2p-1}}\left(\frac{1}{A^{n+k}_{H}(t)}\right)^{\frac{1}{2p-1}}.
\eea
Separating of variables and integrating from $r$ to $R$, we yield
\begin{flalign*}
&\left(\frac{A_{f}(R,x)}{A^{n+k}_{H}(R)}\right)^{\frac{1}{2p-1}} - \left(\frac{A_{f}(r,x)}{A^{n+k}_{H}(r)}\right)^{\frac{1}{2p-1}}& \\
&\quad \leq \left(\frac{n+k-1}{(2p-1)(2p-n-k)}\right)^{\frac{p-1}{2p-1}}\left(||Ric^{\, \mu, H}_{f-}||_{p, \, f}(R)\right)^{\frac{p}{2p-1}}\int^{R}_{r}\frac{(\sqrt{H})^{\frac{n+k-1}{2p-1}}}{sin^{2}(t\sqrt{H})} \, dt.&
\end{flalign*}
\end{pot1}

 Now we prove Theorem \ref{thm2} using Theorem \ref{thm3.1} above. The argument is similar. \\

\begin{pot2}
Using
\bea
\frac{V_{f}(x,r)}{V^{n+k}_{H}(r)} = \frac{\int^{r}_{0}A_{f}(x,t)\, dt}{\int^{r}_{0}A^{n+k}_{H}(t)\, dt},
\eea
we compute
\begin{flalign} \label{eq3.7}
\frac{d}{dr}\left(\frac{V_{f}(x,r)}{V^{n+k}_{H}(r)}\right) = \frac{A_{f}(x,r)\int^{r}_{0}A^{n+k}_{H}(t)\, dt - A^{n+k}_{H}(r)\int^{r}_{0}A_{f}(x,t)\,dt}{\left(V^{n+k}_{H}(r)\right)^{2}}.
\end{flalign}
On the other hand, integrating both sides of (\ref{eq3.4.1}) from $t$ to $r$, we have
\begin{flalign*}
&\frac{A_{f}(x,r)}{A^{n+k}_{H}(r)} - \frac{A_{f}(x,t)}{A^{n+k}_{H}(t)}& \\
&\quad \leq \left((2p-1)^{p}\left(\frac{n+k-1}{2p-n-k}\right)^{p-1}\right)^{\frac{1}{2p-1}}\frac{\left(||Ric^{\, \mu,H}_{f-}||_{p, \, f}(r)\right)^{\frac{p}{2p-1}}}{\left(A^{n+k}_{H}(t)\right)^{1 - \frac{1}{2p-1}}\left(A^{n+k}_{H}(t)\right)^{\frac{1}{2p-1}}}& \\
&\qquad \times \int^{r}_{t}\left(A_{f}(x,s)\right)^{1 - \frac{1}{2p-1}} ds& \\
&\quad \leq \left((2p-1)^{p}\left(\frac{n+k-1}{2p-n-k}\right)^{p-1}\right)^{\frac{1}{2p-1}}\frac{\left(||Ric^{\, \mu,H}_{f-}||_{p, \, f}(r)\right)^{\frac{p}{2p-1}}}{A^{n+k}_{H}(t)}\left(r - t \right)^{\frac{1}{2p-1}} & \\
&\qquad \times \left(V_{f}(x,r)\right)^{1 - \frac{1}{2p-1}}& \\
&\quad \leq \left((2p-1)^{p}\left(\frac{n+k-1}{2p-n-k}\right)^{p-1}\right)^{\frac{1}{2p-1}}\frac{\left(||Ric^{\, \mu,H}_{f-}||_{p, \, f}(r)\right)^{\frac{p}{2p-1}}}{A^{n+k}_{H}(t)} r^{\frac{1}{2p-1}} \left(V_{f}(x,r)\right)^{1 - \frac{1}{2p-1}}.&
\end{flalign*} 
This implies that
\begin{flalign*}
&A_{f}(x,r)\int^{r}_{0}A^{n+k}_{H}(t) \, dt - A^{n+k}_{H}(r)\int^{r}_{0}A_{f}(x,t) \, dt& \\
&\quad \leq \left((2p-1)^{p}\left(\frac{n+k-1}{2p-n-k}\right)^{p-1}\right)^{\frac{1}{2p-1}}r^{\frac{2p}{2p-1}} \left(||Ric^{\, \mu,H}_{f-}||_{p, \, f}(r)\right)^{\frac{p}{2p-1}}A^{n+k}_{H}(r)& \\
&\qquad \times \left(V_{f}(x,r)\right)^{1 - \frac{1}{2p-1}}.&
\end{flalign*}
Inserting this inequality into (\ref{eq3.7}) gives
\begin{flalign*}
&\frac{d}{dr}\left(\frac{V_{f}(x,r)}{V^{n+k}_{H}(r)}\right)& \\
&\quad \leq \left((2p-1)^{p}\left(\frac{n+k-1}{2p-n-k}\right)^{p-1}\right)^{\frac{1}{2p-1}} \frac{\left(||Ric^{\, \mu,H}_{f-}||_{p, \, f}(r)\right)^{\frac{p}{2p-1}}}{\left(V^{n+k}_{H}(r)\right)^{2}}& \\
&\qquad \times A^{n+k}_{H}(r)\left(V_{f}(x,r)\right)^{1 - \frac{1}{2p-1}}r^{\frac{2p}{2p-1}} \\
&\quad \leq \left((2p-1)^{p}\left(\frac{n+k-1}{2p-n-k}\right)^{p-1}\right)^{\frac{1}{2p-1}}\left(||Ric^{\, \mu,H}_{f-}||_{p, \, f}(r)\right)^{\frac{p}{2p-1}}A^{n+k}_{H}(r)& \\
&\qquad \times \left(\frac{V_{f}(x,r)}{V^{n+k}_{H}(r)}\right)^{1 - \frac{1}{2p-1}} \left(\frac{r}{V^{n+k}_{H}(r)}\right)^{\frac{2p}{2p-1}}.&
\end{flalign*}
Separating of variables, we get
\begin{flalign*}
&\frac{d}{dr}\left(\frac{V_{f}(x,r)}{V^{n+k}_{H}(r)}\right)^{\frac{1}{2p-1}}& \\
&\quad \leq \left(\frac{n+k-1}{(2p-1)(2p-n-k)}\right)^{\frac{p-1}{2p-1}}\left(||Ric^{\, \mu,H}_{f-}||_{p, \, f}(r)\right)^{\frac{p}{2p-1}}A^{n+k}_{H}(r)\left(\frac{r}{V^{n+k}_{H}(r)}\right)^{\frac{2p}{2p-1}}.&
\end{flalign*}
Integrating from $r$ to $R$, we obtain
\begin{flalign*}
&\left(\frac{V_{f}(x,R)}{V^{n+k}_{H}(R)}\right)^{\frac{1}{2p-1}} - \left(\frac{V_{f}(x,r)}{V^{n+k}_{H}(r)}\right)^{\frac{1}{2p-1}}& \\
&\quad \leq \left(\frac{n+k-1}{(2p-1)(2p-n-k)}\right)^{\frac{p-1}{2p-1}}\left(||Ric^{\, \mu,H}_{f-}||_{p, \, f}(R)\right)^{\frac{p}{2p-1}}& \\
&\qquad \times \int^{R}_{r}A^{n+k}_{H}(t)\left(\frac{t}{V^{n+k}_{H}(t)}\right)^{\frac{2p}{2p-1}}dt& \\
&\quad \leq \left(\frac{n+k-1}{(2p-1)(2p-n-k)}\right)^{\frac{p-1}{2p-1}}\left(||Ric^{\, \mu,H}_{f-}||_{p, \, f}(R)\right)^{\frac{p}{2p-1}}& \\
&\qquad \times \int^{R}_{0}A^{n+k}_{H}(t)\left(\frac{t}{V^{n+k}_{H}(t)}\right)^{\frac{2p}{2p-1}}dt.& 
\end{flalign*}
Hence, we complete the proof of Theorem \ref{thm2}.
\end{pot2}

 Next, As a corollary we have the volume doubling result.

\begin{corollary} \label{coro3.2}
Let $(M,g,e^{-f}dv)$ be an $n$-dimensional smooth metric measure space. Assume that $\mu \geq \frac{1}{k}$ for some positive constant $k$. For $\beta > 1$ and $p > \frac{n+k}{2}$, there is an $\epsilon = \epsilon(n+k,p,H,R,\beta)$ such that if $\bar{k}(p,H,R) < \epsilon$, then for all $x \in M$ and $0 < r_{1} < r_{2} \leq R$ (assume $R \leq \frac{\pi}{2\sqrt{H}}$ when $H > 0$), we have
\be \label{eq3.8}
\frac{V_{f}(x,r_{2})}{V_{f}(x,r_{1})} \leq \beta \cdot \frac{V^{n+k}_{H}(r_{2})}{V^{n+k}_{H}(r_{1})}.
\ee
\end{corollary}

\begin{poc1}
By Theorem \ref{thm2}, we obtain
\bea
\left(\frac{V_{f}(x,r_{1})}{V^{n+k}_{H}(r_{1})}\right)^{\frac{1}{2p-1}} \geq \left(\frac{V_{f}(x,r_{2})}{V^{n+k}_{H}(r_{2})}\right)^{\frac{1}{2p-1}} - C(n+k,p,H,r_{2})\left(||Ric^{\, \mu,H}_{f-}||_{p, \, f}(r_{2})\right)^{\frac{p}{2p-1}}.
\eea
Multiplying $\left(\frac{V^{n+k}_{H}(r_{1})}{V_{f}(x,r_{2})}\right)^{\frac{1}{2p-1}}$ on the above inequality, then we get
\be \label{eq3.9}
\left(\frac{V_{f}(x,r_{1})}{V_{f}(x,r_{2})}\right)^{\frac{1}{2p-1}} \geq \left(\frac{V^{n+k}_{H}(r_{1})}{V^{n+k}_{H}(r_{2})}\right)^{\frac{1}{2p-1}}\left(1 - \sigma(r_{2})\right),
\ee
where 
\bea
\sigma(r_{2}) &:=& C(n+k,p,H,r_{2})V^{n+k}_{H}(r_{2})^{\frac{1}{2p-1}} \\
&&\times \left(\frac{1}{V_{f}(x,r_{2})}\int_{B(x,r_{2})}|Ric^{\, \mu,H}_{f-}|^{p}A_{f}\, dtd\theta_{n-1}\right)^{\frac{1}{2p-1}}.
\eea
Similarly, we have for $R$
\be \label{eq3.10}
\left(\frac{V_{f}(x,r_{2})}{V_{f}(x,R)}\right)^{\frac{1}{2p-1}} \geq \left(\frac{V^{n+k}_{H}(r_{2})}{V^{n+k}_{H}(R)}\right)^{\frac{1}{2p-1}}\left(1 - \sigma(R)\right),
\ee
where 
\bea
\sigma(R) &:=& C(n+k,p,H,R)V^{n+k}_{H}(R)^{\frac{1}{2p-1}} \\
&&\times \left(\frac{1}{V_{f}(x,R)}\int_{B(x,R)}|Ric^{\, \mu,H}_{f-}|^{p}A_{f}\, dtd\theta_{n-1}\right)^{\frac{1}{2p-1}}. 
\eea
Since $C(n+k,p,H,r)$ is increasing in $r$, we obtain
\be \label{eq3.11}
\sigma(r_{2}) \leq \sigma(R)\left(\frac{V^{n+k}_{H}(r_{2})}{V^{n+k}_{H}(R)}\right)^{\frac{1}{2p-1}}\left(\frac{V_{f}(x,R)}{V_{f}(x,r_{2})}\right)^{\frac{1}{2p-1}}.
\ee
Note that
\bea
\sigma(r_{2}) &\leq& C(n+k,p,H,r_{2})V^{n+k}_{H}(r_{2})^{\frac{1}{2p-1}} \\
&&\times \left(\sup_{x \in M}\left(\frac{1}{V_{f}(x,r_{2})}\int_{B(x,r_{2})}|Ric^{\, \mu,H}_{f-}|^{p}A_{f}\, dtd\theta_{n-1}\right)^{\frac{1}{p}}\right)^{\frac{p}{2p-1}} \\
&=& C(n+k,p,H,r_{2})V^{n+k}_{H}(r_{2})^{\frac{1}{2p-1}}\bar{k}(p,H,r_{2})^{\frac{p}{2p-1}}.
\eea
Similarly, we also have
\bea
\sigma(R) \leq C(n+k,p,H,R)V^{n+k}_{H}(R)^{\frac{1}{2p-1}}\bar{k}(p,H,R)^{\frac{p}{2p-1}}.
\eea
Let $\epsilon(n+k,p,H,R,\beta)^{\frac{p}{2p-1}} = \frac{\left(1 - \left(\frac{1}{\beta}\right)^\frac{1}{2p-1}\right)}{3C(n+k,p,H,R)V^{n+k}_{H}(R)^{\frac{1}{2p-1}}}$. Substituting this on the above inequality, we obtain
\bea
\sigma(R) &\leq& C(n+k,p,H,R)V^{n+k}_{H}(R)^{\frac{1}{2p-1}}\frac{2}{3C(n+k,p,H,R)V^{n+k}_{H}(R)^{\frac{1}{2p-1}}} \\
&=& \frac{2}{3}.
\eea
Putting this into (\ref{eq3.10}) yields
\bea
\left(\frac{V_{f}(x,r_{2})}{V_{f}(x,R)}\right)^{\frac{1}{2p-1}} \geq \frac{1}{3}\left(\frac{V^{n+k}_{H}(r_{2})}{V^{n+k}_{H}(R)}\right)^{\frac{1}{2p-1}}.
\eea
By (\ref{eq3.11}), we get
\bea
\sigma(r_{2}) \leq 3\sigma(R).
\eea
This implies that
\bea
\left(\frac{V_{f}(x,r_{1})}{V_{f}(x,r_{2})}\right)^{\frac{1}{2p-1}} \geq \left(\frac{V^{n+k}_{H}(r_{1})}{V^{n+k}_{H}(r_{2})}\right)^{\frac{1}{2p-1}}\left(1 - 3\sigma(R)\right).
\eea
On the other hand, using
\bea
\epsilon(n+k,p,H,R,\beta)^{\frac{p}{2p-1}} = \frac{\left(1 - \left(\frac{1}{\beta}\right)^\frac{1}{2p-1}\right)}{3C(n+k,p,H,R)V^{n+k}_{H}(R)^{\frac{1}{2p-1}}},
\eea
we have
\bea
\sigma(R) \leq \frac{1 - \left(\frac{1}{\beta}\right)^{\frac{1}{2p-1}}}{3}.
\eea
Hence, we can deduce
\bea
\frac{V_{f}(x,r_{2})}{V_{f}(x,r_{1})} \leq \beta \cdot \frac{V^{n+k}_{H}(r_{2})}{V^{n+k}_{H}(r_{1})},
\eea
where $\beta > 1$.
\end{poc1}

Finally, we study the volume comparison estimate for annulus. The idea of proof is similar to Theorem \ref{thm2}. \\

\begin{pot3}
By (\ref{eq3.7}), we have
\begin{flalign} \label{eq3.12}
\frac{d}{dR}\left(\frac{V_{f}(x,r,R)}{V^{n+k}_{H}(r,R)}\right) = \frac{A_{f}(x,R)\int^{R}_{r}A^{n+k}_{H}(t)\, dt - A^{n+k}_{H}(R)\int^{R}_{r}A_{f}(x,t) \, dt}{\left(V^{n+k}_{H}(r,R)\right)^{2}}.
\end{flalign}
Similarly, we also have
\begin{flalign} \label{eq3.13}
&\frac{d}{dr}\left(\frac{V_{f}(x,r,R_{1})}{V^{n+k}_{H}(r,R_{1})}\right)& \\
&\quad = \frac{-A_{f}(x,r)\int^{R_{1}}_{r}A^{n+k}_{H}(t)\, dt + A^{n+k}_{H}(r)\int^{R_{1}}_{r}A_{f}(x,t) \, dt}{\left(V^{n+k}_{H}(r,R_{1})\right)^{2}}. \nonumber
\end{flalign}
Integrating (\ref{eq3.4.1}) from $t$ to $R$ yields
\begin{flalign*}
&\left(\frac{A_{f}(x,R)}{A^{n+k}_{H}(R)}\right) - \left(\frac{A_{f}(x,t)}{A^{n+k}_{H}(t)}\right)& \\
&\quad \leq \left((2p-1)^{p}\left(\frac{n+k-1}{2p-n-k}\right)^{p-1}\right)^{\frac{1}{2p-1}}\left(||Ric^{\, \mu,H}_{f-}||_{p, \, f}(R)\right)^{\frac{p}{2p-1}}\left(\frac{1}{A^{n+k}_{H}(t)}\right)& \\
&\qquad \times\int^{R}_{t}A_{f}(x,s)^{1 - \frac{1}{2p-1}}ds& \\
&\quad \leq \left((2p-1)^{p}\left(\frac{n+k-1}{2p-n-k}\right)^{p-1}\right)^{\frac{1}{2p-1}}\left(||Ric^{\, \mu,H}_{f-}||_{p, \, f}(R)\right)^{\frac{p}{2p-1}} \left(\frac{1}{A^{n+k}_{H}(t)}\right)&\\
&\qquad \times (R - t)^{\frac{1}{2p-1}}\left(\int^{R}_{t} A_{f}(x,s) \, ds\right)^{1 - \frac{1}{2p-1}}& \\
&\quad \leq \left((2p-1)^{p}\left(\frac{n+k-1}{2p-n-k}\right)^{p-1}\right)^{\frac{1}{2p-1}}\left(||Ric^{\, \mu,H}_{f-}||_{p, \, f}(R)\right)^{\frac{p}{2p-1}} \left(\frac{1}{A^{n+k}_{H}(t)}\right)& \\
&\qquad \times R^{\frac{1}{2p-1}}V_{f}(x,t,R)^{1 - \frac{1}{2p-1}}.&
\end{flalign*}
Multiplying $A^{n+k}_{H}(R)A^{n+k}_{H}(t)$ on the above inequality, then we obtain
\begin{flalign*}
&A_{f}(x,R)A^{n+k}_{H}(t) - A_{f}(x,t)A^{n+k}_{H}(R)& \\
&\quad \leq \left((2p-1)^{p}\left(\frac{n+k-1}{2p-n-k}\right)^{p-1}\right)^{\frac{1}{2p-1}}\left(||Ric^{\, \mu,H}_{f-}||_{p, \, f}(R)\right)^{\frac{p}{2p-1}}R^{\frac{1}{2p-1}}& \\
&\qquad \times A^{n+k}_{H}(R) V_{f}(x,t,R)^{1 - \frac{1}{2p-1}}.&
\end{flalign*}
Integrating this inequality from $r$ to $R$, we get
\begin{flalign} \label{eq3.14}
&A_{f}(x,R)\int^{R}_{r}A^{n+k}_{H}(t) \, dt - A^{n+k}_{H}(R)\int^{R}_{r}A_{f}(x,t) \, dt& \\
&\quad \leq \left((2p-1)^{p}\left(\frac{n+k-1}{2p-n-k}\right)^{p-1}\right)^{\frac{1}{2p-1}}\left(||Ric^{\, \mu,H}_{f-}||_{p, \, f}(R)\right)^{\frac{p}{2p-1}}  \nonumber \\
&\qquad \times (R-r)R^{\frac{1}{2p-1}} A^{n+k}_{H}(R) V_{f}(x,r,R)^{1 - \frac{1}{2p-1}} \nonumber \\
&\quad \leq \left((2p-1)^{p}\left(\frac{n+k-1}{2p-n-k}\right)^{p-1}\right)^{\frac{1}{2p-1}}\left(||Ric^{\, \mu,H}_{f-}||_{p, \, f}(R)\right)^{\frac{p}{2p-1}} \nonumber \\
&\qquad \times R^{\frac{2p}{2p-1}} A^{n+k}_{H}(R) V_{f}(x,r,R)^{1 - \frac{1}{2p-1}}. \nonumber
\end{flalign}
On the other hand, similar to the argument, we also get
\begin{flalign} \label{eq3.15}
&A^{n+k}_{H}(r)\int^{R_{1}}_{r}A_{f}(x,t) \, dt - A_{f}(x,r)\int^{R_{1}}_{r}A^{n+k}_{H}(t) \, dt& \\
&\quad \leq \left((2p-1)^{p}\left(\frac{n+k-1}{2p-n-k}\right)^{p-1}\right)^{\frac{1}{2p-1}}\left(||Ric^{\, \mu,H}_{f-}||_{p, \, f}(R_{1})\right)^{\frac{p}{2p-1}} \nonumber \\
&\qquad \times R_{1}^{\frac{2p}{2p-1}}A^{n+k}_{H}(R_{1}) V_{f}(x,r,R_{1})^{1 - \frac{1}{2p-1}}. \nonumber
\end{flalign}
Substituting (\ref{eq3.14}) into (\ref{eq3.12}) gives
\bea
\frac{d}{dR}\left(\frac{V_{f}(x,r,R)}{V^{n+k}_{H}(r,R)}\right)^{\frac{1}{2p-1}} &\leq& \left(\frac{n+k-1}{(2p-1)(2p-n-k)}\right)^{\frac{p-1}{2p-1}}\left(||Ric^{\, \mu,H}_{f-}||_{p, \, f}(R)\right)^{\frac{p}{2p-1}} \\
&&\times \left(\frac{R}{V^{n+k}_{H}(r,R)}\right)^{\frac{2p}{2p-1}}A^{n+k}_{H}(R).
\eea
Separating of variables, integrating this inequality from $R_{1}$ to $R_{2}$ and changing the variable $r$ to $r_{2}$, we have
\begin{flalign} \label{eq3.16}
&\left(\frac{V_{f}(x,r_{2},R_{2})}{V^{n+k}_{H}(r_{2},R_{2})}\right)^{\frac{1}{2p-1}} - \left(\frac{V_{f}(x,r_{2},R_{1})}{V^{n+k}_{H}(r_{2},R_{1})}\right)^{\frac{1}{2p-1}} \\
&\quad \leq \left(\frac{n+k-1}{(2p-1)(2p-n-k)}\right)^{\frac{p-1}{2p-1}}\int^{R_{2}}_{R_{1}}A^{n+k}_{H}(t)\left(\frac{t}{V^{n+k}_{H}(r_{2},t)}\right)^{\frac{2p}{2p-1}}dt \nonumber \\
&\qquad \times \left(||Ric^{\, \mu,H}_{f-}||_{p, \, f}(R_{2})\right)^{\frac{p}{2p-1}}. \nonumber
\end{flalign}
On the other hand, combining (\ref{eq3.13}) and (\ref{eq3.15}), we obtain
\bea
\frac{d}{dr}\left(\frac{V_{f}(x,r,R_{1})}{V^{n+k}_{H}(r,R_{1})}\right)^{\frac{1}{2p-1}} &\leq& \left(\frac{n+k-1}{(2p-1)(2p-n-k)}\right)^{\frac{p-1}{2p-1}}\left(||Ric^{\, \mu,H}_{f-}||_{p, \, f}(R_{1})\right)^{\frac{p}{2p-1}} \\
&& \times A^{n+k}_{H}(R_{1}) \left(\frac{R_{1}}{V^{n+k}_{H}(r,R_{1})}\right)^{\frac{2p}{2p-1}}.
\eea
Integrating the above inequality from $r_{1}$ to $r_{2}$, we get
\begin{flalign} \label{eq3.17}
&\left(\frac{V_{f}(x,r_{2},R_{1})}{V^{n+k}_{H}(r_{2},R_{1})}\right)^{\frac{1}{2p-1}} - \left(\frac{V_{f}(x,r_{1},R_{1})}{V^{n+k}_{H}(r_{1},R_{1})}\right)^{\frac{1}{2p-1}}& \\
&\quad \leq \left(\frac{n+k-1}{(2p-1)(2p-n-k)}\right)^{\frac{p-1}{2p-1}}\left(||Ric^{\, \mu,H}_{f-}||_{p, \, f}(R_{1})\right)^{\frac{p}{2p-1}} A^{n+k}_{H}(R_{1}) \nonumber \\
&\qquad \times \int^{r_{2}}_{r_{1}}\left(\frac{R_{1}}{V^{n+k}_{H}(t,R_{1})}\right)^{\frac{2p}{2p-1}}dt. \nonumber
\end{flalign}
Adding (\ref{eq3.16}) and (\ref{eq3.17}) gives 
\begin{flalign*}
&\left(\frac{V_{f}(x,r_{2},R_{2})}{V^{n+k}_{H}(r_{2},R_{2})}\right)^{\frac{1}{2p-1}} - \left(\frac{V_{f}(x,r_{1},R_{1})}{V^{n+k}_{H}(r_{1},R_{1})}\right)^{\frac{1}{2p-1}}& \\
&\quad \leq \left(\frac{n+k-1}{(2p-1)(2p-n-k)}\right)^{\frac{p-1}{2p-1}}\left(||Ric^{\, \mu,H}_{f-}||_{p, \, f}(R_{1})\right)^{\frac{p}{2p-1}}& \\
&\qquad \times \left(\int^{R_{2}}_{R_{1}}A^{n+k}_{H}(t)\left(\frac{t}{V^{n+k}_{H}(r_{2},t)}\right)^{\frac{2p}{2p-1}}dt+ \int^{r_{2}}_{r_{1}}A^{n+k}_{H}(R_{1})\left(\frac{R_{1}}{V^{n+k}_{H}(t,R_{1})}\right)^{\frac{2p}{2p-1}}dt\right) &
\end{flalign*}
for $0 \leq r_{1} \leq r_{2} \leq R_{1} \leq R_{2}$.
\end{pot3}

\section{Diameter estimate}
 In this section, using mean curvature comparison estimate and volume comparison estimates, we prove the global diameter estimate. \\

 Let $p_{1}$ and $p_{2}$ are two points in $M$ and $x_{0}$ be a midpoint between $p_{1}$ and $p_{2}$. Consider the excess function
\bea
e(x) = d(p_{1},x) + d(p_{2},x) - d(p_{1},p_{2}).
\eea
Note that $e(x) \geq 0$ in $M$ and $e(x) \leq 2r$ on $B(x_{0},r)$ by the triangle inequality. Using this fact, we want to prove our result by contradiction. That is, we will show that the excess function $e$ is negative on $B(x_{0},r)$.

 From the mean curvature comparison estimate, by using a suitably lare comparison sphere we may choose any large $D$ enough so that if $d(p_{1},p_{2}) > D$, then we have $\Delta_{f} \, e \leq -K + \psi_{1}$ on $B(x_{0},r)$, where $K$ is a large positive constant to be determined, and $\psi_{1}$ denotes an error term controlled by $C_{1}(n+k,p,H,r)\cdot\bar{k}(p,H,r)$.

 Let $\Omega_{j}$ be a sequence of smooth star-shaped domains which converges to $B(x_{0},r) - \mbox{\rm Cut}(x_{0})$ and $u_{i}$ be a sequence of smooth functions such that $|u_{i} - e| < i^{-1}$, $|\nabla u_{i}| \leq 2 + i^{-1}$, and $\Delta_{f} \, u_{i} \leq \Delta_{f} \, e + i^{-1}$ on $B(x_{0},r)$ (see \cite{CO}). Set $h = d^{2}(x_{0},\cdot) - r^{2}$, we have that $h$ is a negative and smooth function on $\Omega_{j}$. By Green's theorem, we obtain
\bea
\int_{\Omega_{j}}\left(\Delta_{f} \, u_{i}\right)h - \int_{\Omega_{j}}\left(\Delta_{f} \, h\right)u_{i} = \int_{\partial\Omega_{j}}h\left(u_{i} \, \nu\right) - \int_{\partial\Omega_{j}}u_{i}\left(h \, \nu\right),
\eea
where $\nu$ is the ourward unit normal to $\Omega_{j}$. Note that
\bea
\Delta_{f} \, u_{i} \leq -K + \psi_{1} + i^{-1}
\eea
and
\bea
\left(\Delta_{f} \, h\right)u_{i} &\leq& \left(e + i^{-1}\right)\left(2d\Delta_{f} \, d + 2\right) \\
&\leq& 3r\left(2(n+k) + \psi_{2}\right),
\eea
where $\psi_{2}$ is another error term controlled by $C_{2}(n+k,p,H,r)\cdot\bar{k}(p,H,r)$. Thus, we have
\bea
\int_{\Omega_{j}}\left(-K + \psi_{1} + i^{-1}\right)h - 3r\int_{\Omega_{j}}\left(2(n+k) + \psi_{2}\right) &\leq& \int_{\partial\Omega_{j}}(i^{-1} - e)(\nu \, h) \\
&& - \int_{\partial\Omega_{j}}h(2 + i^{-1}).
\eea
Since $u_{i} \rightarrow e$ when $i \rightarrow \infty$, by the dominated convergence theorem, the above inequality becomes
\bea
\int_{\Omega_{j}}(-K + \psi_{1})h - 3r\int_{\Omega_{j}}(2n + 2k + \psi_{2}) \leq -2\int_{\partial\Omega_{j}}h - \int_{\partial\Omega_{j}}e(\nu \, h).
\eea
Also, we compute that
\bea
\int_{B(x_{0},r)}(-K + \psi_{1})h &=& -\int_{B(x_{0},r)}Kh + \int_{B(x_{0},r)}\psi_{1}h \\
&\geq& -\int_{B(x_{0},\frac{r}{2})}Kh + \int_{B(x_{0},r)}\psi_{1}h \\
&\geq& \int_{B(x_{0},\frac{r}{2})}K\frac{3}{r}r^{2} - \int_{B(x_{0},r)}\psi_{1}r^{2} \\
&=& \frac{3}{4}r^{2}KV_{f}\left(x_{0},\frac{r}{2}\right) - \int_{B(x_{0},r)}\psi_{1}r^{2}.
\eea
So we get
\begin{flalign*}
&\int_{B(x_{0},r)}(-K + \psi_{1})h - 3r\int_{B(x_{0},r)}(2n + 2k + \psi_{2})& \\
&\quad \geq \frac{3}{4}r^{2}KV_{f}\left(x_{0},\frac{r}{2}\right) - 6r(n+k)V_{f}\left(x_{0},\frac{r}{2}\right)- \int_{B(x_{0},r)}r^{2}\psi_{1} + 3r\psi_{2}.&
\end{flalign*}
By the volume comparison estimate, we obtain
\bea
V_{f}\left(x_{0},\frac{r}{2}\right) \geq 2^{-1}\frac{V^{n+k}_{H}(\frac{r}{2})}{V^{n+k}_{H}(r)}V_{f}(x_{0},r).
\eea
Moreover, if $\bar{k}(p,H,r)$ is small, we also have
\bea
\int_{B(x_{0},r)}(r^{2}\psi_{1} + 3r\psi_{2}) \leq (r^{2} + nr)V_{f}(x_{0},r).
\eea
Hence, we get
\begin{flalign*}
&\int_{B(x_{0},r)}(-K + \psi_{1})h - 3r\int_{B(x_{0},r)}(2n + 2k + \psi_{2})& \\
&\quad \geq r^{2} V_{f}(x_{0},r) \left(\frac{3}{8}\frac{V^{n+k}_{H}(\frac{r}{2})}{V^{n+k}_{H}(r)}K -1 -7r^{-1}(n+k)\right).&
\end{flalign*}
Thus for $K > \frac{8}{3}\left(\frac{V^{n+k}_{H}(r)}{V^{n+k}_{H}(\frac{r}{2})}\right)(7r^{-1}(n+k) + 1)$, the above inequality is positive. This implies that
\bea
-2\int_{\partial\Omega_{j}}h - \int_{\partial\Omega_{j}}e(\nu \, h) > 0
\eea
as $j \rightarrow \infty$. However the first integral term goes to $0$ as $j \rightarrow \infty$, while the second ingeral term : $\nu \, h \geq 0$ on $\partial\Omega_{j}$  for all $j$. This implies that $e$ must be nagative on $B(x_{0},r)$, which is a contradiction. So $d(p_{1},p_{2}) < D$ for some $D$.
\hfill\fbox{}\par
\end{document}